\documentstyle[12pt]{article}

\newcommand{\Two}[4]{\left(\begin{array}{cc}#1&#2\\#3&#4\end{array}\right)}
\newcommand{\Three}[9]{\left(\begin{array}{ccc}#1&#2&#3\\
#4&#5&#6\\#7&#8&#9\end{array}\right)}
\newcommand{\EN}{\end{equation}}
\renewcommand{\Re}{\mbox{\rm Re}}
\renewcommand{\det}{\mbox{det }}
\def\les{\leq}

\def\GL{\mbox{GL}}
\def\G{\mbox{G}}
\def\N{\mbox{N}}
\def\Z{\mbox{Z}}
\def\A{{\bf A}}
\def\C{{\bf C}}
\def\R{{\bf R}}

\begin{document}

\title{On poles of twisted tensor L-functions}

\author{Yuval Z. Flicker and Dmitrii Zinoviev
\thanks {
The first named author wishes to express his very deep gratitude to
Professor Toshio Oshima for inviting him to visit the University of
Tokyo in March 1995, and to him and Professor Takayuki Oda for their
warm hospitality.\hfil\break
Department of Mathematics, The Ohio State University,
231 W. 18th Ave., Columbus, OH 43210-1174. }
\smallskip\\}

\date{}
\maketitle

\vspace{-0.5in}
\begin{abstract}

It is shown that the only possible pole of the twisted tensor
$L$-functions in Re$(s)\geq 1$ is located at $s=1$ for {\em all}
quadratic extensions of global fields.

\end{abstract}

\vskip0pt plus1in
\penalty-1000
\vskip0pt minus1in

\begin{center}{\bf  0. Introduction.}\end{center}

Let $E$ be a quadratic separable field extension of a global 
field $F$. Denote by $\A_E$, $\A_F$ the corresponding rings of adeles.
Put $\G_n$ for $\GL_n$ and $\Z_n$ for its center. Then $\Z_n(\A_E)$
is the group $\A_E^\times$ of ideles of $\A_E$. 
Fix a cuspidal representation
$\pi$ of the adele group $\G_n(\A_E)$. Without lost of generality,
we may assume that the central character of $\pi$ is trivial
on the split component of $\A_E^\times$. This is the
multiplicative group $\R^{\times}$ of the field of real 
numbers embedded in $\A_E^{\times}$ via $x\mapsto (x,...,x,1,...)$
($x$ in the archimedean, $1$ in the finite components).
Let $S$ be a finite
set of places of $F$ (depending on $\pi$), including
the places where $E/F$ ramify, and the archimedean places, such
that for each place $v'$ of $E$ above a place $v$ outside
$S$ the component $\pi_{v'}$ of $\pi$ is unramified.
Following [1], let $r$ be the twisted tensor representation of 
$\widehat{\G}=[\GL(n,\C)\times\GL(n,\C)]\times\mbox{Gal}(E/F)$ on 
$\C^n\otimes\C^n$. It acts by $r((a,b))(x\otimes y)=ax\otimes by$
and $r(\sigma)(x\otimes y)=y\otimes x$ ($\sigma\in\mbox{Gal}(E/F),\,
\sigma\neq 1$).
Let $q_v$ be the cardinality of the residue field $R_v/
\mbox{\boldmath  $\pi$}_vR_v$
of the ring $R_v$ of integers in $F_v$. We define the twisted
tensor $L$-function to be the Euler product
$$
L(s,\,r(\pi),\,S)=\prod_{v\in\!\!\!/ S}\det[1-q_v^{-s}r(t_v)]^{-1}.
$$
\indent
The representation $\pi$ is called {\em distinguished} if its central 
character is trivial on $\A_F^{\times}$ and there is an automorphic 
form $\phi\in\pi$ in $L^2(\G_n(E)\Z_n(\A_F)\backslash\G_n(\A_E))$, such that 
$\int\phi(g)dg\neq 0$. The integral is taken over the closed
subspace 
$\G_n(\!F)\Z_n(\A_F\!)\!\backslash\!\G_n(\A_F)$ of $\G_n(\!E)\Z_n(\!\A_F)
\!\backslash\!\G_n(\A_E\!)$.
\newline\indent
The following theorem is proven in [1, p. 309] for a quadratic extension
$E/F$ of global fields, such that  each archimedean place of $F$ 
splits in $E$. We prove it for any quadratic extension of global 
fields, i.e. also for number fields with completions $E_v/F_v=\C/\R$.

{\bf Theorem.} {\em The product $L(s,\,r(\pi),\,S)$ converges absolutely,
uniformly in compact subsets, in some right half-plane. It has analytic
continuation as a meromorphic function to the right half plane
$\Re(s)>1-\epsilon$, for some small $\epsilon>0$.
The only possible pole of
$L(s,\,r(\pi),\,S)$ in $\Re(s)> 1-\epsilon$ is simple, located at $s=1$.
The function $L(s,\,r(\pi),\,S)$ has a pole at $s=1$ if and only if 
$\pi$ is distinguished. }

{\em Proof.} The proof of this theorem is the same as that of the
Theorem of [1, \S 4], pp. 309-310. On lines 14 and 18 of page 310
of [1], we use the proposition below. It holds in the non-split 
archimedean case too. Hence the restriction put in [1] on the extension
$E/F$ can be removed.

For the functional equation satisfied by $L(s,\,r(\pi),\,S)$, see [1].
For the local $L$-factors at all non-archimedean places of $F$, see 
[2]. The non-vanishing of this $L$-function on the edge $\Re(s)=1$ 
of the critical strip has been shown by Shahidi [6].
Twisted tensor L-functions are used in the study (see Kon-no [5]) of 
the residual spectrum of unitary groups.
\begin{center}{\bf  1. Local computations.}\end{center}

From now on, we consider the local case only.
Let $E/F$ be a quadratic extension of local fields. 
Thus in the archimedean case $E/F\,=\,\C/\R$.
Denote by $x\mapsto\bar{x}$ the non-trivial automorphism
of $E$ over $F$. Let $\iota\neq 0$ be an element of $E$, such that
$\bar{\iota}=-\iota$.
Put $\G_n$ for $\GL_n$. The groups of $F$ and $E$-points are denoted
by $\G_n(F)$ and $\G_n(E)$. Denote by $\N_n$ the unipotent radical
of the upper triangular subgroup of $\G_n$, and by A$_n$ the diagonal
subgroup. Let $\psi_0$ be a non trivial additive character of $F$.
For example, if $F=\R$ then $\psi_0(x)=e^{2\pi ix}$.
Let $\psi$ be the (non-trivial) character $\psi(z)=
\psi_0((z-\bar{z})/\iota)$ of $E$. It is trivial on $F$.
For $u\in\N_n(E)$, set $\theta(u)=\psi(\sum_{i=1}^{n-1}u_{i,i+1})$.
\newline\indent
Fix an irreducible admissible representation $\pi$ of $\G_n(E)$ on
a complex vector space $V$. The representation $\pi$ is called
{\em generic} if there exists a non-zero linear form $\lambda$ on
V, such that $\lambda(\pi(u)v)=\theta(u)\lambda(v)$ for all $v$ in $V$
and $u$ in $\N_n(E)$. The dimension of the space of such $\lambda$
is bounded by one. Let $W(\pi;\,\theta)$ be the space of functions $W$
on $\G_n(E)$ of the form $W(g)=\lambda(\pi(g)v)$, where $v\in V$.
We have $W(ug)=\theta(u)W(g)$ ($g\in\G_n(E)$, $u\in\N_n(E)$).
Denote by $W_0(\pi;\,\theta)$ those functions in $W(\pi;\,\theta)$
whose corresponding vectors $v$ are in the space of K-finite 
vectors, where K$=$K$_n(E)$ is the standard maximal compact 
subgroup of $\G_n(E)$. 

For $\Phi\in S(F^n)$, define the integral
$$
\Psi(s,\Phi,W)=\int_{N_n(F)\backslash G_n(F)}W(g)\Phi(\epsilon_n g)
|\det g|^sdg, $$ 
where $\epsilon_n=(0,0,...,0,1)$ is a row vector
of size $n$.

{\bf Proposition.}
 (i) {\em  There exists some small constant $\epsilon$, $\epsilon >0$, 
such that the integral 
$\Psi(s,\Phi,W)$ converges absolutely, uniformly in compact subsets, 
for $\Re(s) > 1-\epsilon$; }  
\newline\indent  
(ii) {\em  There exists $W$ in $W_0(\pi;\,\theta)$ and $\Phi$ in
$S(F^n)$, such that $\Psi(s,\Phi,W)\neq 0$.}

{\em Proof.} When $E/F$ is an extension of non-archimedean local
fields, (i) and (ii) are treated in the Proposition of [1], \S 4, p. 308.
We prove (i) in general, including the case $(E,\,F)=(\C,\,\R)$, following
Jacquet and Shalika [3], pp. 204-206. 
\newline\indent
Using the Iwasawa decomposition $\G_n(F)= 
\N_n(F)\mbox{A}_n(F)\mbox{K}_n(F)$,  
and the associated measure decomposition, we 
need to show the convergence of the integral
$$
\int_{A_n(F)K_n(F)}|W(ak)||\det a|^{s}\delta_{n,F}^{-1}
 (a)dadk.
$$
Here $a=$diag$(a_1,\,a_2,\,...\,,a_{n-1},\,1)$. Recall that
$$
\delta_{n,F}(a)=\delta_{n-1,F}(a)|\det a|=|\det a|
\prod_{1\le i<j\le n-1}\frac{|a_i|}{|a_j|},
$$
and (see e.g. [1], p. 307) that $\delta_{n,E}(a)=\delta_{n,F}^2(a)$.
\newline\indent
By Proposition 3 of Jacquet and Shalika [3, \S 4] there is a 
finite set ${\bf X}$ of finite
functions in $n-1$ variables such that $|W(ak)|$ is bounded by a finite
sum of expressions of the form
$$ C_\chi\delta_{n-1,E}^{1/2}(a)\Phi\left(\frac{a_1}{a_2},\,
\frac{a_2}{a_3},\,...,\,a_{n-1}\right). $$
Here $C_\chi$ is the absolute value of some element of ${\bf X}$ and
$\Phi\ge 0$ is in $S(F^{n-1})$. Thus, it suffices to show that the
integral obtained by replacing $W$ by this estimate is convergent.
Using that
$$ 
\delta_{n-1,E}^{1/2}(a)\delta_{n,F}^{-1}(a)=
\delta_{n-1,F}(a)\delta_{n,F}^{-1}(a)=|\det a|^{-1},
$$
we arrive at the finite sum of integrals
$$ \int C_\chi\Phi\left(\frac{a_1}{a_2},\,\frac{a_2}{a_3},\,...,
   \,a_{n-1}\right)|\det a|^{s-1}da.
$$
The change of variables
$ a_1=t_1...t_{n-1},\,a_2=t_2...t_{n-1},\,...,\,a_{n-1}=t_{n-1},$
has the Jacobian $t_2t_3^2...t_{n-2}^{n-3}$. We obtain a sum of
expressions of the form
$$ \int C_\chi\Phi(t_1,\,t_2,\,...,\,t_{n-1})\prod_{j=2}^{n-2}t_j^{j-1}
   \prod_{j=1}^{n-1}t_j^{j(s-1)}dt.
$$
Again, by Proposition 3 of Jacquet and Shalika [3, \S 4] the set 
${\bf X}$ is such that any $\chi$ in it 
is the product of (1) a polynomial in the logarithms of the
absolute values of the variables, and (2) a character of the form
$$
\chi_1(t_1)\chi_2(t_2)...\chi_{n-1}(t_{n-1}),
$$
with $\Re(\chi_i)>0$, for each $i$. It follows that the above 
integral converges uniformly in compact subsets of $\Re(s)>1-\epsilon$, 
for some small $\epsilon>0$. This completes the proof of (i).
\newline\indent
  For (ii) we will follow the proof of Proposition 7.3 of Jacquet and
Shalika [3]. Assume that $\Psi(s,\Phi,W)=0$ for all choices of $W$ 
in $W_0(\pi;\,\theta)$ and $\Phi$ in $S(F^n)$. We will show that
$W(e)=0$ for all $W$, a contradiction which will imply (ii) of the lemma.
Since $\Phi$ is arbitrary, it follows that for all $W$ we have
$$\int_{N_{n-1}(F)\backslash G_{n-1}(F)}W\left[\Two{g}{0}{0}{1}\right]
|\det g|^sdg=0.$$
Define
$$ I_k(W)=\int_{N_k(F)\backslash G_k(F)}W
\left[\Two{g}{0}{0}{1_{n-k}}\right]|\det g|^sdg.$$
We claim that $I_k(W)$ is zero for all $W$ and all $k$ with $0\les k \les n-1$.
The lemma would then follow, since $W(e)=I_0(W)$. We will show 
this claim by descending induction on $k$.
We have just seen that $I_{n-1}(W)=0$.
So fix $k\les n-1$ with $I_k(W)=0$ for all $W$. We proceed to show that 
$I_{k-1}(W)=0$ for all $W$.
\newline\indent
We apply the fact that $I_k(W)=0$ to the function $W_\Phi$ defined by
$$
W_\Phi(g)=\int_{F^k}W\left[g\Three{1_k}{\iota u}{0}{0}{1}{0}{0}{0}
{1_{n-k-1}}\right]\Phi(u)du. 
$$
Here $u$ is a column of size $k$, $\Phi\in S(F^k)$ and $W\in W_0(\pi;\theta)$.
Proposition 2.4 of Jacquet and Shalika [4; II], p. 784, and the remark 
following it (top of p.
786), assure us that this function is in the space $W_0(\pi;\,\theta)$.
\newline\indent
Note that
$$W_\Phi\left[\Two{g}{0}{0}{1_{n-k}}\right]=
W\left[\Two{g}{0}{0}{1_{n-k}}\right]\widehat{\Phi}(\epsilon_k g),$$
where $ \widehat{\Phi}(y)=\int_{F^k} \Phi(u)\psi_0(y\cdot u)du$ 
denotes the Fourier transform of $\Phi\in S(F^k)$.
Indeed
$$ \widehat{\Phi}(\epsilon_k g)=\int_{F^k}\Phi(u)\psi_0(\epsilon_k 
g\cdot u)du  = \int_{F^k}\Phi(u)\psi_0\left(\sum_{j=1}^k g_{kj}u_j
\right)du.$$
Further, since
$$\Three{g}{0}{0}{0}{1}{0}{0}{0}{1_{n-k-1}}
  \Three{1_k}{\iota u}{0}{0}{1}{0}{0}{0}{1_{n-k-1}}
\, = \,  \Three{1_k}{\iota gu}{0}{0}{1}{0}{0}{0}{1_{n-k-1}}
  \Three{g}{0}{0}{0}{1}{0}{0}{0}{1_{n-k-1}},$$
we have
\begin{eqnarray*}
\lefteqn{W\left[\Three{1_k}{\iota gu}{0}{0}{1}{0}{0}{0}{1_{n-k-1}}
  \Three{g}{0}{0}{0}{1}{0}{0}{0}{1_{n-k-1}}\right]} \\
                                                   \\
 &=& \theta\left(\Three{1_k}{\iota gu}{0}{0}{1}{0}{0}{0}{1_{n-k-1}}\right)
  W\left[\Three{g}{0}{0}{0}{1}{0}{0}{0}{1_{n-k-1}}\right]  \\
                                                   \\
 &=& \psi\left(\sum_{j=1}^k \iota g_{kj}u_j\right)
     W\left[\Three{g}{0}{0}{0}{1}{0}{0}{0}{1_{n-k-1}}\right]
 =  \psi_0\left(\sum_{j=1}^k g_{kj}u_j\right)
     W\left[\Three{g}{0}{0}{0}{1}{0}{0}{0}{1_{n-k-1}}\right].
\end{eqnarray*}
Now substituting $W_\Phi$ for $W$ in $I_k(W)=0$, we obtain
$$\int_{N_k(F)\backslash G_k(F)}W\left[\Two{g}{0}{0}{1_{n-k}}\right]
\widehat{\Phi}(\epsilon_k g)|\det g|^s dg=0$$
for all $\Phi\in S(F^k)$ and all 
$W\in W_0(\pi;\,\theta)$. In this integral $\widehat{\Phi}$ 
can be replaced by any element of $S(F^k)$.
Hence $I_{k-1}(W)=0$ for all $W$ and we are done.

\begin{center}{\bf References}\end{center}

[1] Y. Flicker, ``Twisted Tensors and Euler Products'', {\em
    Bull. Soc. Math. France} 116, 1988, 295--313.

[2] Y. Flicker, ``On the local twisted tensor $L$-function'', 
    Appendix to: ``On zeroes of the twisted tensor $L$-function'', 
    {\em Math. Ann.} 297, 1993, 199--219.

[3] H. Jacquet and J. Shalika, ``Exterior square $L$-functions'', in {\em
    Automorphic Forms, Shimura Varieties, and $L$-functions}, 1990, 
    143--226 (ed. L. Clozel, S. Milne).

[4] H. Jacquet and J. Shalika, ``On Euler products and the classification
    of automorphic representations'', {\em Amer. J. Math.} 103, 1981, 
    I, 449--558; II, 777--815.

[5] T. Kon-no, ``The residual spectrum of $U(2,2)$, preprint.

[6] F. Shahidi, ``On certain $L$-functions'', {\em Amer. J. Math.} 103,
    1981, 297--355.

\end{document}